\documentclass[10pt]{amsart}
\usepackage{amssymb}
\usepackage{epsfig}
\usepackage{url}
\theoremstyle{plain}

\newtheorem{thm}{Theorem}[section]

\newtheorem{rem}[thm]{Remark}

\numberwithin{equation}{section}

\def\bbb{\mathbb}
\renewcommand{\phi}{\varphi}

\newcommand{\Q}{\bbb{Q}}

\begin{document}
\title{On arithmetic progressions on genus two curves }
\author{Maciej Ulas}
\thanks{The author is a scholar of a project which is
co-financed from the European Social Fund and Polish national
budget within the Integrated Regional Operational Programme.}

\date{}

\begin{abstract}
We study arithmetic progression in the $x$-coordinate of rational
points on genus two curves. As we know, there are two models for
the curve $C$ of genus two: $C:\;y^2=f_{5}(x)$ or
$C:\;y^2=f_{6}(x)$, where $f_{5},\;f_{6}\in\Q[x]$,
$\operatorname{deg}f_{5}=5,\;\operatorname{deg}f_{6}=6$ and the
polynomials $f_{5},\;f_{6}$ do not have multiple roots. First we
prove that there exists an infinite family of curves of the form
$y^2=f(x)$, where $f\in\Q[x]$ and $\operatorname{deg}f=5$ each
containing 11 points in arithmetic progression. We also present an
example of $F\in\Q[x]$ with $\operatorname{deg}F=5$ such that on
the curve $y^2=F(x)$ twelve points lie in arithmetic progression.
Next, we show that there exist infinitely many curves of the form
$y^2=g(x)$ where $g\in\Q[x]$ and $\operatorname{deg}g=6$, each
containing 16 points in arithmetic progression. Moreover, we
present two examples of curves in this form with 18 points in
arithmetic progression.
\end{abstract}

\maketitle
\bigskip

\section{Introduction}

Let $f\in\Q[X]$ be a polynomial without multiple roots and let us
consider the curve $C:\;y^2=f(x)$. We say that rational points
$P_{i}=(x_{i},\;y_{i})$ for $i=1,\;2,\;\ldots,\;n$ are in {\it
arithmetic progression} on the curve $C$, if rational numbers
$x_{i}$ are in arithmetic progression for $i=1,\;2,\;\ldots,\;n$.
A positive integer $n$ will be called {\it the length of
arithmetic progression} on the curve $C$. A natural question
arises here: How long can arithmetic progression on the curve
$y^2=f(x)$ with fixed degree of $f$ be? Through the whole paper by
a point we mean a rational one.

In case of polynomials of degree one, this question is equivalent
to the question about the number of squares which form an
arithmetic progression.

It is not difficult to show that there exists an infinite family
$\mathcal{A}_{1}$ of polynomials of degree one, with the property
that for each $f\in\mathcal{A}_{1}$ there are 3 points in
arithmetic progression on the curve $y^2=f(x)$ (of genus 0). It
turns out, however, which was already proved by Fermat, that it is
impossible to construct arithmetic progression composed of four
squares.

In paper \cite{A} Allison has shown that there exists an infinite
family $\mathcal{A}_{2}$ of polynomials of degree two, such that
for each $f\in\mathcal{A}_{2}$ on the curve $y^2=f(x)$ (of genus
0) eight points lie in arithmetic progression.

In the case of polynomials of degree three, Bremner in \cite{Br}
has constructed an infinite family $\mathcal{A}_{3}$ with such a
property that for every $f\in\mathcal{A}_{3}$ on the curve
$y^2=f(x)$ (of genus 1) 8 points lie in arithmetic progression. A
similar result with the use of other methods was obtained by
Campbell in \cite{Ca}.

In the case of polynomials of degree four, in \cite{Ul} we have
constructed an infinite family $\mathcal{A}_{4}$ with such
property that for every $f\in\mathcal{A}_{4}$ there are 12 points
in arithmetic progression on the curve $y^2=f(x)$ (of genus 1).

It is worth noting that MacLeod in \cite{ML} has constructed
polynomials $F_{i},\;(i=1,2,3,4)$ of degree four such that on each
curve $y^2=F_{i}(x)$ there are 14 points in arithmetic
progression.

In all above cases, each of the families
$\mathcal{A}_{2},\;\mathcal{A}_{3},\;\mathcal{A}_{4}$ is
parametrized by rational points on some elliptic curve of positive
rank.

It is reasonable to define the following quantities

\begin{center}
$m(d):=\operatorname{max}\{k:$ there exists a polynomial
$g\in\Q[x]$ of degree $\operatorname{deg}g=d$ such that
  \\\hskip 2.2cm on the curve $y^2=g(x)$ there are $k$ points in arithmetic progression$\}$,
\end{center}

\begin{center}
$M(d):=\operatorname{max}\{k:$ there exists an infinite family
$\mathcal{A}_{d}$ of polynomials of degree $d$,\\ \hskip 2cm such
that for every $g\in\mathcal{A}_{d}$  there are
\\\hskip -0.8cm $k$ points in arithmetic progression on the curve $y^2=g(x)$ $\}$.
 \end{center}

We have obvious inequality $m(d)\geq M(d)$. The above results can
be grouped in the following manner:

\bigskip

\begin{center}
            \begin{tabular}{|c|c|c|c|c|}
              \hline
              $d$ & $1$ & $2$ & $3$ & $4$ \\
              \hline
              $m(d)$ & $3$ & $\geq 8$ & $\geq 8$ & $\geq 14$ \\
              $M(d)$ & $3$ & $\geq 8$ & $\geq 8$ & $\geq 12$ \\
              \hline
            \end{tabular}
         \end{center}
                 \begin{center}{\sc Table 1}\end{center}
\bigskip

 In this paper we will concentrate on quantities $m(d)$ and
$M(d)$ for $d=5,\;6$. Let us note that it corresponds to the
construction of arithmetic progressions on hyperelliptic curves of
genus 2. In case of $d=5$ we show that $m(5)\geq 12$ and $M(5)\geq
11$. When $d=6$, we first show that there exists a polynomial
$G(t,x)\in\Q(t)[x]$ such that  there are 14 points in arithmetic
progression on the curve $y^2=G(t,x)$. Using another approach we
prove that $m(6)\geq 18$ and $M(6)\geq 16$.

\bigskip

\section{Case of $d=5$}

Using a method similar to that used by Campbell in \cite{Ca} we
will show the following

\bigskip

\begin{thm}\label{thm1}
         There exist polynomials $F_{i}(t,x)\in\Q(t)[x],\;(i=1,2)$ of degree
         $\operatorname{deg}_{x}F_{i}=5$ such that on the
         curve $y^2=F_{i}(t,x)$ eleven $\Q(t)$-rational points lie in
         arithmetic progression.
\end{thm}
\begin{proof}
           Let $u$ be a variable and let us consider a polynomial
           \begin{equation*}
              g(u,x)=(x-u)^2\prod_{i=1}^{10}(x-i).
\end{equation*}
As we know, there is exactly one pair of polynomials
$h,\;f\in\Q(u)[x]$ such that
$\operatorname{deg}_{x}h=6,\;\operatorname{deg}_{x}f=5$ and
\begin{equation*}
              g(u,x)=h(u,x)^2-\frac{25}{1048576}f(u,x).
\end{equation*}
In our case the polynomial $f$ is in the form of
\begin{equation*}
                 f(u,x)=a_{5}x^5+a_{4}x^4+a_{3}x^3+a_{2}x^2+a_{1}x+a_{0},
  \end{equation*}
 where
\begin{align*}
             &a_{0}=5695244944u^2-12894461800u+263250625,\\
             &a_{1}=-8(533634200u^2-873304794u-1611807725), \\
             &a_{2}=32(37257330u^2-8034400u-396302603),\\
             &a_{3}=-3520(41536u^2+145226u-1285845),\\
             &a_{4}=3520(1888u^2+31152u-193477),\\
             &a_{5}=-(6645760u-36551680).
 \end{align*}
 If $u$ is rational and $u\neq 11/2$ then the polynomial $f(u,x)$ is
 without multiple roots. Thus, we see that in this case there are 10 points in arithmetic progression on the
 curve $y^2=f(u,x)$.

 Let us now consider curve $Q_{1}$ with the equation
 $p_{1}^2=f(u_{1},11)$. This is a quadric with rational point
 $(11,\;16225)$. Using the standard method we have a parametrization
 of the curve $Q_{1}$ given by:

 \begin{align*}
    &p_{1}(t)=\frac{16225(5695244944t^2-794728t+1)}{5695244944t^2-1},\\
    &u_{1}(t)=\frac{11(4523021144t^2+2950t-1)}{5695244944t^2-1}.
   \end{align*}

   If we now define $F_{1}(t,x)=f(u_{1}(t),x)$, where $u_{1}(t)$
   is as above, then on the curve
   \begin{equation*}
   C_{1}:\;y^2=F_{1}(t,x)
   \end{equation*}
    there are 11 points in arithmetic progression.

   We can similarly parametrize the quadric $Q_{2}$ given by the
   equation $p_{2}^2=f(u_{2},0)$ with rational point $(0,16225)$.
   In this case the parametrization takes the form
    \begin{align*}
    &p_{2}(t)=\frac{16225(5695244944t^2+794728t+1)}{5695244944t^2-1},\\
    &u_{2}(t)=\frac{32450t(397364t+1)}{5695244944t^2-1}.
   \end{align*}
   If we now define $F_{2}(t,x)=f(u_{2}(t),x)$, where $u_{2}(t)$
   is as above, then on the curve
   \begin{equation*}
   C_{2}:\;y^2=F_{2}(t,x)
   \end{equation*}
   there are 11 points in arithmetic progression.
   \end{proof}

\bigskip

   Finding of such a rational $t$ that on the curve $C_{1}$ there are twelve points
   in arithmetic progression requires finding rational points on
   the curve
   \begin{align*}
       y^2=&-543542815457978537904123051776t^4-5858532530788995918150400t^3\\
          &+611541611111856733408t^2-633115875308400t-30556659591,
   \end{align*}
   or on the curve
   \begin{align*}
       y^2=&10452723797211797241575306232064t^4+16006824835104105921670400t^3\\
       &-4781502606421467214112t^2-3647410080111600t+547548809049.
   \end{align*}
   Then, points with $x$-coordinates in $\{1,\;2,\;\ldots,\;12\}$
   (respectively with $x$-coordinates in
   $\{0,\;1,\;\ldots,\;11\}$) will be in arithmetic progression on
   the curve $C_{1}$. In the case of the curve $C_{2}$ we obtain the
   same curves. It is easy to see that the above curves have
   $\Q_{p}$-rational points for every $p$, but unfortunately we
   did not manage to find a rational point on any of the above curves.
   It seems that finding a rational point on any of the curves (or showing
   that such points do not exist) may be a difficult task.

\begin{rem}
{\rm The statement of Theorem \ref{thm1} can also be obtained
using the following reasoning. Let us consider the polynomial
\begin{equation*}
f(t,x)=(x-t)\prod_{i=1}^{11}(x-i).
\end{equation*}
Then there exist polynomials $p,\;F\in\Q(t)[x]$, such that
$\operatorname{deg}_{x}p=6,\operatorname{deg}_{x}F=5$ and
\begin{equation*}
f(t,x)=p(t,x)^2-F(t,x)
\end{equation*}
Then the curve $C:\;y^2=F(t,x)$ contains 11 points in arithmetic
progression. Unfortunately, in this case the polynomials $F(t,0)$
and $F(t,12)$ are irreducible of degree 6, and each of the curves
$y^2=F(t,0),\;y^2=F(t,12)$ contain only finitely many rational
points. Therefore, the curve $C$ cannot be used to construct an
infinite family of curves with the required property.}
\end{rem}

\bigskip

  The following example found with the use of computer shows that
  $m(5)\geq 12$. Consider the curve
  \begin{equation*}
                  C:\;y^2=\;12x^5-322x^4+3208x^3-14438x^2+27980x-16079.
    \end{equation*}
    We have the following points in arithmetic progression on the
    curve $C$:
    \begin{align*}
    \{&(1, 19),\;(2, 55),\;(3, 37),\;(4, 1),\;(5, 11),\;(6, 31),\;(7, 35),\\
               &(8, 23),\;(9, 29),\;(10, 89),\;(11, 181),\;(12, 305)\}.
    \end{align*}

\bigskip

\section{Case of $d=6$}

Let us begin with the following

\bigskip

\begin{thm}\label{thm2}
There exists a polynomial $H(t,x)\in\Q(t)[x]$ of degree
$\operatorname{deg}_{x}H=6$, such that fourteen $\Q(t)$-rational
points lie in arithmetic progression on the curve $y^2=H(t,x)$.
\end{thm}
\begin{proof}
Let $t$ be a variable and let us consider a polynomial
\begin{equation*}
                h(t,x)=(x^2-15x+4t)\prod_{i=1}^{14}(x-i).
 \end{equation*}
 Then there is exactly one pair of polynomials
$g,\;H\in\Q(t)[x]$ such that
$\operatorname{deg}_{x}8=6,\;\operatorname{deg}_{x}H=6$ and
\begin{equation*}
              h(t,x)=g(t,x)^2-H(t,x).
\end{equation*}
In our case the polynomial $H$ is in the form
\begin{equation*}
                H(t,x)=a_{3}(x(15-x))^3+a_{2}(x(15-x))^2+a_{1}(x(15-x))+a_{0},
  \end{equation*}
  where
  \begin{align*}
             a_{0}=&\;46228440064 - 37262033920t + 10620980224t^2 - 1420209280t^3\\
                   &+106891216t^4 - 4876960t^5 + 144760t^6 - 2800t^7 + 25t^8,\\
             a_{1}=&\;4(-790888960 + 642389312t - 177526160t^2 + 21803240t^3\\
                   & -1364540t^4 + 42826t^5 - 630t^6 + 5t^7),\\
             a_{2}=&\;2(35503616 - 29056640t + 7910592t^2 - 929040u^3+52318t^4\\
                    &-1260t^5+7t^6),\\
             a_{3}=&\;2(-261120+215008t-58040t^2+6636t^3-350t^4+7t^5).
 \end{align*}
 Therefore, we see that on the curve
 \begin{equation*}
                 C:\;y^2=H(t,x)
 \end{equation*}
  fourteen points lie in arithmetic progression. These points
 are of the form $P_{i}=(i,\;g(t,i))$ for
 $i=1,\;2,\;\ldots,\;14.$
\end{proof}

\bigskip

The first part of the proof of Theorem \ref{thm2} suggests
considering polynomials which are invariant with respect to the
change of variables $x\rightarrow 15-x$. Let us, therefore,
consider the polynomial
\begin{equation}\label{R1}
              f(x)=b_{3}(x(x-15))^3+b_{2}(x(x-15))^2+b_{1}x(x-15)+b_{0},
    \end{equation}
    where
    \begin{align*}
               &b_{0}=(6p^2-22q^2+27r^2-11s^2)/47520,\\
               &b_{1}=(159p^2-517q^2+567r^2-209s^2)/11880,\\
               &b_{2}=(5496p^2-14872q^2+14337r^2-4961s^2)/11880,\\
               &b_{3}=(156p^2-308q^2+273r^2-91s^2)/30.
    \end{align*}
    For $f$ defined in this way we have
    \begin{equation*}
      f(1)=f(14)=p^2,\;f(2)=f(13)=q^2,\;f(3)=f(12)=r^2,\;f(4)=f(11)=s^2.
    \end{equation*}
    Therefore, we see that in order to obtain on the curve $y^2=f(x)$ an arithmetic
    progression of the length 14, it is necessary to investigate a
    system of equations
    \begin{equation}\label{R3}
     \begin{cases}
             f(5)=(-14p^2+77q^2-162r^2+154s^2)/55=u^2 \\
             f(6)=(-21p^2+110q^2-210r^2+154s^2)/33=v^2\\
             f(7)=(-60p^2+308q^2-567r^2+385s^2)/66=w^2
     \end{cases}
     \end{equation}

     Using a substitution
     $(p,\;q,\;r,\;s,\;u)=(a+u,\;b+u,\;c+u,\;d+u,\;u)$ we obtain
     a parametrization of solutions of the first equation of system
     (\ref{R3})
     \begin{align*}
               (p,\;q,\;r,\;s,\;u)=(&14a^2-154ab+77b^2+324ac-162c^2-308ad+154d^2,\\
                                    &14a^2-28ab+77b^2-324bc+162c^2+308bd-154d^2,\\
                                    &-14a^2+77b^2+28ac-154bc+162c^2-308cd+154d^2,\\
                                    &14a^2-77b^2+162c^2-28ad+154bd-324cd+154d^2,\\
                                    &-14a^2+77b^2-162c^2+154d^2).
              \end{align*}
Now let us set
\begin{equation}\label{R4}
(a,\;b,\;c,\;d)= (946A,\;946,\;11(15A+71),\;441A+505).
\end{equation}
For $a,\;b,\;c,\;d$ defined in this way, we get a parametric
solution of the system (\ref{R3}) given by
\begin{align*}
(p,\;q,\;r,\;s,\;u,\;v,\;w)=(&181144A^2+85170A-42585,59140A^2-118280A-164589,\\
                             &17230A^2-112505A-128454,52874A^2+102845A+68010,\\
                             &59140A^2+122004A+42585,43984A^2+104070A+75675,\\
                             &15790A^2+107955A+99984)
\end{align*}
   For $p,\;q,\;r,\;s$ defined above, the coefficients of
   the polynomial $g_{A}(x)=36f(x)$ are
   \begin{align*}
       b_{0}=&\;36(128941675300A^4+235814377620A^3+34730973441A^2\\
                  &\;-216866857320A-132565503600),\\
       b_{1}=&\;4(A-1)(254A+219)(354070194A^2+848446325A+620203644),\\
       b_{2}=&\;(A-1)(254A+219)(35708622A^2+96399845A+73722213),\\
       b_{3}=&\;4(A-1)(254A+219)(72474A^2+210275A+164709).
       \end{align*}
      For $A\in\Q\setminus S$, where $S=\{-240/233,-219/254,1,475/2\},$ the polynomial $g_{A}$
       does not have multiple roots. From this we can conclude
       that for $A\in\Q\setminus S$ on the curve
       \begin{equation*}
                      C_{A}:\;y^2=g_{A}(x)
       \end{equation*}
       fourteen points lie in arithmetic progression. Now, it
       is an easy task to prove the following

\bigskip

\begin{thm}\label{thm3}
         There exist infinitely many $A\in\Q$ such that on the
         curve $C_{A}:\;y^2=g_{A}(x)$ there are 16 points in
         arithmetic progression.
\end{thm}
\begin{proof}
           Let us set $x=0$ and consider the curve
           \begin{equation*}
                          C:\;y^2=b_{0}(A).
           \end{equation*}
It is easy to see that on $C$ we have rational point
$P=(1,1342374)$. As we know, the curve of the form $y^2=f_{4}(x)$,
where $\operatorname{deg}f_{4}=4$, with rational point is
birationally equivalent to an elliptic curve with Weierstrass'
equation \cite{Mor}. Using APECS program \cite{Co} we obtain that
$C$ is birational with the curve
\begin{equation*}
 E:y^2+xy+y=x^3-x^2+21015110653x+1214962664541571.
\end{equation*}
For the curve $E$ we have
\begin{equation*}
                \operatorname{Tors}E(\mathbb{Q})=\{\mathcal{O},\;(-51365,25682)\},
 \end{equation*}
 and again using APECS we obtain that free part of $E(\Q)$ is
 generated by
 \begin{align*}
   &G_{1}=(-45989,\; -12274606),\;G_{2}=(751451,\;-664705966),\\
   &G_{3}=(-17669,\; -28941646),\;G_{4}=(24913585/256,\; 264676595567/4096).
   \end{align*}
   As an immediate consequence, we get that there are infinitely
   many rational points on the curve $C$ and all but finitely many
   define the curve $C_{A}:\;y^2=g_{A}(x)$ with 16 points in
   arithmetic progression.
\end{proof}

\bigskip

To show that $m(6)\geq 18$ we have taken the polynomial of the
form
\begin{equation}\label{R5}
 h(x)=c_{3}(x(x-19))^3+c_{2}(x(x-19))^2+c_{1}x(x-19)+c_{0}.
\end{equation}
With help of computer we found the following numbers
$c_{0},\;c_{1},\;c_{2},\;c_{3}$ such that the polynomial
(\ref{R5}) has values which are squares of integers for
$x=1,\;2,\;,\ldots,\;18$:
\bigskip
\begin{center}
\begin{tabular}{|c|c|c|c|}
  \hline
  $c_{0}$ & $c_{1}$ & $c_{2}$ & $c_{3}$ \\
  \hline
  358043904 & 18892800 & 321792 & 1664 \\
  864002304 & 37085184 & 524544 & 2432 \\
  \hline
\end{tabular}
\end{center}
            \begin{center}{\sc Table 2}\end{center}

 \vskip 1cm

\noindent {\bf Acknowledgments.} I would like to thank the
anonymous referee for his valuable comments and Professor K. Rusek
for remarks improving the presentation.

\vskip 0.5cm

 \hskip 4.5cm       Maciej Ulas

 \hskip 4.5cm       Jagiellonian University

 \hskip 4.5cm       Institute of Mathematics

 \hskip 4.5cm       Reymonta 4

 \hskip 4.5cm       30 - 059 Krak\'{o}w, Poland

 \hskip 4.5cm      e-mail:\;{\tt Maciej.Ulas@im.uj.edu.pl}
\end{document}